
%
%
%
%
%
%
\magnification=\magstephalf      
%
%
\vsize=7.5truein                 
\hsize=5.2truein                 
\newskip\stdskip                 
\stdskip=6pt plus3pt minus3pt    
\medskipamount=\stdskip          
\parindent=0pt                   
\parskip=\stdskip                
\abovedisplayskip=\stdskip       
\belowdisplayskip=\stdskip       
\mathsurround=0.75pt             
\overfullrule=0pt                
%
%
\def\ppar{\par\goodbreak\vskip 8pt plus 4pt minus 4pt}     
%
%
\def\stdspace{\hskip 0.75em plus 0.15em\ignorespaces}
\let\qua\stdspace 
%
%
%
%
%
%
%
\def\hexnumber#1{\ifcase#1 0\or 1\or 2\or 3\or 4\or 5\or 6\or 7\or 8\or
 9\or A\or B\or C\or D\or E\or F\fi}
%
%
\font\thirtnmsa=msam10 scaled 1315    
\font\tenmsa=msam10          \font\ninemsa=msam9
\font\sevenmsa=msam7         \font\sixmsa=msam6
\font\fivemsa=msam5
%
%
\newfam\msafam                  \textfont\msafam=\tenmsa
\scriptfont\msafam=\sevenmsa    \scriptscriptfont\msafam=\fivemsa
\edef\hexa{\hexnumber\msafam}        
\def\msa{\fam\msafam\tenmsa}         
%
%
\font\thirtnmsb=msbm10 scaled 1315   
\font\tenmsb=msbm10      \font\ninemsb=msbm9
\font\sevenmsb=msbm7     \font\sixmsb=msbm6
\font\fivemsb=msbm5
%
\newfam\msbfam                   \textfont\msbfam=\tenmsb       
\scriptfont\msbfam=\sevenmsb     \scriptscriptfont\msbfam=\fivemsb
\edef\hexb{\hexnumber\msbfam}    
\def\msb{\fam\msbfam\tenmsb}     
%
%
\font\thirtneufm=eufm10 scaled 1315   
\font\teneufm=eufm10                 \font\nineeufm=eufm9
\font\seveneufm=eufm7                \font\sixeufm=eufm6
\font\fiveeufm=eufm5
%
\newfam\eufmfam                    \textfont\eufmfam=\teneufm
\scriptfont\eufmfam=\seveneufm     \scriptscriptfont\eufmfam=\fiveeufm
\edef\hexf{\hexnumber\eufmfam}      
\def\frak{\fam\eufmfam\teneufm}     
%
%
%
\font\thirtnrm=cmr10 scaled 1315    
\font\ninerm=cmr9                   \font\sixrm=cmr6   
%
\font\thirtni=cmmi10 scaled 1315    
\font\ninei=cmmi9                   \font\sixi=cmmi6  
%
\font\thirtnsy=cmsy10 scaled 1315   
\font\ninesy=cmsy9                  \font\sixsy=cmsy6  
%
\font\thirtnbf=cmbx10 scaled 1315   
\font\ninebf=cmbx9                  \font\sixbf=cmbx6  
%
%
\font\thirtnex=cmex10 scaled 1315   
\font\nineex=cmex9                  
%
%
\font\thirtnit=cmti10 scaled 1315  
\font\nineit=cmti9                  
%
\font\thirtnsl=cmsl10 scaled 1315  
\font\ninesl=cmsl9                  
%
\font\thirtntt=cmtt10 scaled 1315  
\font\ninett=cmtt9                  
%
%
%
%
\def\small{%
%
%
\textfont0=\ninerm \scriptfont0=\sixrm \scriptscriptfont0=\fiverm
\def\rm{\fam0\ninerm}
%
%
\textfont1=\ninei \scriptfont1=\sixi \scriptscriptfont1=\fivei
%
%
\textfont2=\ninesy \scriptfont2=\sixsy \scriptscriptfont2=\fivesy
%
%
\textfont3=\nineex \scriptfont3=\nineex \scriptscriptfont3=\nineex
%
%
\textfont\bffam=\ninebf \scriptfont\bffam=\sixbf
\scriptscriptfont\bffam=\fivebf \def\bf{\fam\bffam\ninebf}%
%
%
\textfont\itfam=\nineit \def\it{\fam\itfam\nineit}%
\textfont\slfam=\ninesl \def\sl{\fam\slfam\ninesl}%
\textfont\ttfam=\ninett \def\tt{\fam\ttfam\ninett}%
%
%
%
\textfont\msafam=\ninemsa \scriptfont\msafam=\sixmsa
\scriptscriptfont\msafam=\fivemsa \def\msa{\fam\msafam\ninemsa}%
%
%
\textfont\msbfam=\ninemsb \scriptfont\msbfam=\sixmsb
\scriptscriptfont\msbfam=\fivemsb \def\msb{\fam\msbfam\ninemsb}%
%
%
\textfont\eufmfam=\nineeufm  \scriptfont\eufmfam=\sixeufm
\scriptscriptfont\eufmfam=\fiveeufm \def\frak{\fam\eufmfam\nineeufm}%
%
%
%
\normalbaselineskip=11pt%
\setbox\strutbox=\hbox{\vrule height8pt depth3pt width0pt}%
%
%
\normalbaselines\rm
%
%
\stdskip=4pt plus2pt minus2pt    
\medskipamount=\stdskip          
\parskip=\stdskip                
\abovedisplayskip=\stdskip       
\belowdisplayskip=\stdskip       
\def\ppar{\par\goodbreak\vskip 6pt plus 3pt minus 3pt}%
%
%
\def\section##1{\global\advance\sectionnumber by 1
\vskip-\lastskip\penalty-800\vskip 20pt plus10pt minus5pt 
\egroup{\bf\number\sectionnumber\quad##1}\bgroup\small         
\vskip 6pt plus3pt minus3pt
\nobreak\resultnumber=1}
}    
%
\def\beginsmall{\bgroup\small}
\let\endsmall\egroup
%
%
%
%
\def\large{%
\textfont0=\thirtnrm \scriptfont0=\ninerm \scriptscriptfont0=\sevenrm
\def\rm{\fam0\thirtnrm}%
\textfont1=\thirtni \scriptfont1=\ninei \scriptscriptfont1=\seveni
\textfont2=\thirtnsy \scriptfont2=\ninesy \scriptscriptfont2=\sevensy
\textfont3=\thirtnex \scriptfont3=\thirtnex \scriptscriptfont3=\thirtnex
\textfont\bffam=\thirtnbf \scriptfont\bffam=\ninebf
\scriptscriptfont\bffam=\sevenbf \def\bf{\fam\bffam\thirtnbf}%
\textfont\itfam=\thirtnit \def\it{\fam\itfam\thirtnit}%
\textfont\slfam=\thirtnsl \def\sl{\fam\slfam\thirtnsl}%
\textfont\ttfam=\thirtntt \def\tt{\fam\ttfam\thirtntt}%
\textfont\msafam=\thirtnmsa \scriptfont\msafam=\ninemsa
\scriptscriptfont\msafam=\sevenmsa \def\msa{\fam\msafam\thirtnmsa}%
\textfont\msbfam=\thirtnmsb \scriptfont\msbfam=\ninemsb
\scriptscriptfont\msbfam=\sevenmsb \def\msb{\fam\msbfam\thirtnmsb}%
\textfont\eufmfam=\thirtneufm  \scriptfont\eufmfam=\nineeufm
\scriptscriptfont\eufmfam=\seveneufm \def\frak{\fam\eufmfam\teneufm}%
\normalbaselineskip=16pt%
\setbox\strutbox=\hbox{\vrule height11.5pt depth4.5pt width0pt}%
\normalbaselines\rm}%
\let\Large\large   
%
\def\Bbb#1{{\msb#1}}

%

%
\mathchardef\plussquare="0\hexa01
\mathchardef\nge="3\hexb0B
\mathchardef\maltesecross="0\hexa7A
\mathchardef\del="0\hexf01
%
%
%
%
\font\sc=cmcsc10
%
%
%
%
\def\sqr#1#2{{\vcenter{\vbox{\hrule  height.#2truept
	\hbox{\vrule width.#2truept height#1truept 
	\kern#1truept \vrule width.#2truept}
	\hrule height.#2truept}}}}
\def\sq{\sqr55}    
%
%
%
%
\newcount\sectionnumber            
\newcount\resultnumber             
\sectionnumber=0\resultnumber=1    
%
%
%
\def\section#1{\global\advance\sectionnumber by 1
\xdef\nextkey{\number\sectionnumber}
\vskip-\lastskip\penalty-800\vskip 20pt plus10pt minus5pt 
{\large\bf\number\sectionnumber\quad#1}         
\vskip 8pt plus4pt minus4pt
\nobreak\resultnumber=1}      
%
%
%
%
%
         
%
%
%
%

%
\def\proc#1{\xdef\nextkey{\number\sectionnumber.\number\resultnumber}%
\vskip-\lastskip\ppar\bf%
\noindent#1\ \number\sectionnumber.\number\resultnumber
\stdspace\sl\global\advance\resultnumber by 1\ignorespaces}
\def\endproc{\rm\ppar} 
%
%
\def\prf{\vskip-\lastskip\ppar\noindent{\bf Proof}%
\stdspace\rm}                            
\def\endprf{\unskip\stdspace\hbox{}
\hfill$\sq$\par\medskip}                 
\def\proof#1{\vskip-\lastskip\ppar\noindent{\bf#1}%
\stdspace\rm\ignorespaces}        
%
%
%
%
%
%
%
%
\def\proclaim#1{\vskip-\lastskip\ppar\bf%
\noindent#1\stdspace\sl\ignorespaces} 

%
%
%
%
\def\rk#1{\vskip-\lastskip\ppar{\bf #1}\stdspace\ignorespaces}                
\def\endrk{\par\medskip}
%
%
%
%
%
%
\def\label{\xdef\nextkey{\number\sectionnumber.\number\resultnumber}%
\number\sectionnumber.\number\resultnumber
\global\advance\resultnumber by 1}
%
%
%
%
%
%
%
%
%
%
%
%
%
%
%
%
\newcount\refnumber              
\refnumber=1                     
\long\def\reflist#1\endreflist{%
\long\def\thereflist{#1}{\def\refkey##1##2\par{\xdef##1{\number\refnumber}%
\global\advance\refnumber by 1}%
\def\key##1##2\par{\expandafter\xdef%
\csname##1\endcsname{\number\refnumber}%
\global\advance\refnumber by 1}#1\par}}
\long\def\references{%
\penalty-800\vskip-\lastskip\vskip 15pt plus10pt minus5pt 
{\large\bf References}\ppar 
{\leftskip=25pt\frenchspacing    
\small\parskip=3pt plus2pt       
\def\refkey##1##2\par{\noindent  
\llap{[##1]\stdspace}\ignorespaces##2\par}         
\def\key##1##2\par{\noindent  
\llap{[\ref{##1}]\stdspace}\ignorespaces##2\par}  
\def\,{\thinspace}\thereflist\par}}
%
%
%
\newcount\footnotenumber         
\footnotenumber=1                
\def\fnote#1{\xdef\nextkey{\number\footnotenumber}%
{\small\ifnum\footnotenumber>9\parindent=14pt%
\else\parindent=10pt\fi\footnote{$^{\number\footnotenumber}$}%
{\hglue-5pt#1}\global\advance\footnotenumber by 1}}
%
%
%
%
%
%
%
\newcount\figurenumber          
\figurenumber=1                 
\def\caption#1{\xdef\nextkey{\number\figurenumber}%
\cl{\small Figure \number\figurenumber: #1}%
\global\advance\figurenumber by 1}
\def\figurelabel{\xdef\nextkey{\number\figurenumber}%
\cl{\small Figure \number\figurenumber}%
\global\advance\figurenumber by 1}
\long\def\figure#1\endfigure{{\xdef\nextkey{\number\figurenumber}%
\let\captiontext\relax\def\caption##1{\xdef\captiontext{##1}}%
\midinsert\cl{\ignorespaces#1\unskip\unskip\unskip\unskip}\vglue6pt\cl{\small 
Figure \number\figurenumber\ifx\captiontext\relax\else: \captiontext
\fi}\endinsert\global\advance\figurenumber by 1}}
%
%
%
%
%
%
%
\def\nextkey{??}   
%
\def\key#1{\expandafter\xdef\csname #1\endcsname{\nextkey}}
\def\ref#1{\expandafter\ifx\csname #1\endcsname\relax
\immediate\write16{Reference {#1} undefined}??\else
\csname #1\endcsname\fi}
%
%
%
%
%
%
%
\newread\gtinfile
\newwrite\gtreffile
\def\useforwardrefs{
\openin\gtinfile\jobname.ref
\ifeof\gtinfile
\closein\gtinfile
\immediate\write16{No file \jobname.ref}
\else
\closein\gtinfile
\input \jobname.ref
\fi
\immediate\openout\gtreffile \jobname.ref
%
%
\def\key##1{{\def\\{\noexpand}%
\expandafter\xdef\csname ##1\endcsname{\nextkey}%
\immediate\write\gtreffile{\\\expandafter\\\def\\\csname ##1\\\endcsname%
{\nextkey}}}}
%
%
\long\def\reflist##1\endreflist{%
\long\def\thereflist{##1}{\def\refkey####1####2\par{\xdef####1{%
\number\refnumber}{\def\\{\noexpand}\immediate\write\gtreffile
{\\\def\\####1{\number\refnumber}}}\global\advance\refnumber by 1}%
\def\key####1####2\par{\expandafter\xdef%
\csname####1\endcsname{\number\refnumber}%
{\def\\{\noexpand}\immediate\write\gtreffile
{\\\expandafter\\\def\\\csname ####1\\\endcsname{\number\refnumber}}}
\global\advance\refnumber by 1}##1\par}}
\long\def\biblio##1\endbiblio{\reflist##1\endreflist\references}%
%
%
\def\numkey##1{{\def\\{\noexpand}%
\xdef##1{\number\sectionnumber.\number\resultnumber}
\immediate\write\gtreffile{\\\def\\##1%
{\number\sectionnumber.\number\resultnumber}}}}
\def\seckey##1{{\def\\{\noexpand}\xdef##1{\number\sectionnumber}
\immediate\write\gtreffile{\\\def\\##1{\number\sectionnumber}}}}
\def\figkey##1{\xdef##1{\number\figurenumber}%
{\def\\{\noexpand}\immediate\write\gtreffile%
{\\\def\\##1{\number\figurenumber}}}
\number\figurenumber\global\advance\figurenumber by 1}
}   
%
%
%
%
\def\figkey#1{\xdef#1{\number\figurenumber}%
\number\figurenumber\global\advance\figurenumber by 1}
\def\fig#1#2\endfig{%
\midinsert\cl{#2}\vglue6pt\cl{\small Figure #1}\endinsert}
\def\newfig{\number\figurenumber\global\advance\figurenumber by 1}
\def\numkey#1{\xdef#1{\number\sectionnumber.\number\resultnumber}}
\def\seckey#1{\xdef#1{\number\sectionnumber}}
%
%
%
%
%
%
%
%
%
\def\verb{\catcode`\"=\active}       
\def\brev{\catcode`\"=12}            
\brev                                
\verb                                
{\obeyspaces\gdef {\ }}              
{\catcode`\`=\active\gdef`{\relax\lq}}
\def"{%
\begingroup\baselineskip=12pt\def\par{\leavevmode\endgraf}%
\tt\obeylines\obeyspaces\parskip=0pt\parindent=0pt%
\catcode`\$=12\catcode`\&=12\catcode`\^=12\catcode`\#=12%
\catcode`\_=12\catcode`\~=12%
\catcode`\{=12\catcode`\}=12\catcode`\%=12\catcode`\\=12%
\catcode`\`=\active\let"\endgroup}
\brev      
%
%
%
%
%
%
\def\items{\par\leftskip = 25pt}           
\def\enditems{\par\leftskip = 0pt}         
\def\item#1{\par\leavevmode\llap{#1\stdspace}%
\ignorespaces}                             
%
%

%
%
\def\np{\vfil\eject}         
\def\nl{\hfil\break}         
\def\cl{\centerline}         
\def\agt{{\mathsurround=0pt\it$\cal A\mskip-.7mu$lgebraic \&\ 
$\cal G\mskip-2mu$eometric $\cal T\!\!$opology}}  
%
%
%

%
%
%
%
%
\def\title#1{\def\thetitle{#1}}

\def\author#1{\edef\previousauthors{\theauthors}
 \ifx\theauthors\relax\def\theauthors{#1}\else
 \def\theauthors{\previousauthors\par#1}\fi}

%
\def\address#1{\edef\previousaddresses{\theaddress}
 \ifx\theaddress\relax\def\theaddress{#1}\else
 \def\theaddress{\previousaddresses\par\vskip 2pt\par#1}\fi}
\def\secondaddress#1{\edef\previousaddresses{\theaddress}
 \ifx\theaddress\relax\def\theaddress{#1}\else
 \def\theaddress{\previousaddresses\par{\rm and}\par#1}\fi}   

\def\email#1{\edef\previousemails{\theemail}
 \ifx\theemail\relax\def\theemail{#1}\else
 \def\theemail{\previousemails\hskip 0.75em\relax#1}\fi}
\def\secondemail#1{\edef\previousemails{\theemail}
 \ifx\theemail\relax\def\theemail{#1}\else
 \def\theemail{\previousemails\hskip 0.75em{\rm and}\hskip 0.75em
 \relax#1}\fi}
\def\url#1{\edef\previousurls{\theurl}
 \ifx\theurl\relax\def\theurl{#1}\else
 \def\theurl{\previousurls\hskip 0.75em\relax#1}\fi}
\def\secondurl#1{\edef\previousurls{\theurl}
 \ifx\theurl\relax\def\theurl{#1}\else
 \def\theurl{\previousurls\hskip 0.75em{\rm and}\hskip 0.75em
 \relax#1}\fi}
\long\def\abstract#1\endabstract{\long\def\theabstract{#1}}
\def\primaryclass#1{\def\theprimaryclass{#1}}
\def\secondaryclass#1{\def\thesecondaryclass{#1}}
\def\keywords#1{\def\thekeywords{#1}}
%
%
\let\\\par\let\thetitle\relax\let\theshorttitle\relax
\let\theauthors\relax\let\theshortauthors\relax
\let\theaddress\relax\let\theshortaddress\relax
\let\theemail\relax\let\theurl\relax
\let\theabstract\relax\let\theprimaryclass\relax
\let\thesecondaryclass\relax\let\thekeywords\relax
%
%
%
%
\long\def\maketitlepage{    

\vglue 0.2truein   

%
{\parskip=0pt\leftskip 0pt plus 1fil\def\\{\par\smallskip}{\large
\bf\thetitle}\par\medskip}   

\vglue 0.15truein 

%
{\parskip=0pt\leftskip 0pt plus 1fil\def\\{\par}{\sc\theauthors}
\par\medskip}%
 
\vglue 0.1truein 

%
{\small\parskip=0pt
{\leftskip 0pt plus 1fil\def\\{\par}{\sl\theaddress}\par}
\ifx\theemail\relax\else  
\vglue 5pt \def\\{\stdspace{\rm and}\stdspace} 
\cl{Email:\stdspace\tt\theemail}\fi
\ifx\theurl\relax\else    
\vglue 5pt \def\\{\stdspace{\rm and}\stdspace} 
\cl{URL:\stdspace\tt\theurl}\fi\par}

\vglue 7pt 

{\bf Abstract}

\vglue 5pt

\theabstract

\vglue 7pt 

{\bf AMS Classification numbers}\quad Primary:\quad \theprimaryclass\par

Secondary:\quad \thesecondaryclass

\vglue 5pt 

{\bf Keywords:}\quad \thekeywords

\np  

}    
%
%
\long\def\makeshorttitle{    


%
{\parskip=0pt\leftskip 0pt plus 1fil\def\\{\par\smallskip}{\large
\bf\thetitle}\par\medskip}   

\vglue 0.05truein 

%
{\parskip=0pt\leftskip 0pt plus 1fil\def\\{\par}{\sc\theauthors}
\par\medskip}%
 
\vglue 0.03truein 

%
{\small\parskip=0pt
{\leftskip 0pt plus 1fil\def\\{\par}{\sl\ifx\theshortaddress\relax
\theaddress\else\theshortaddress\fi}\par}
\ifx\theemail\relax\else  
\vglue 5pt \def\\{\stdspace{\rm and}\stdspace} 
\cl{Email:\stdspace\tt\theemail}\fi
\ifx\theurl\relax\else    
\vglue 5pt \def\\{\stdspace{\rm and}\stdspace} 
\cl{URL:\stdspace\tt\theurl}\fi\par}

\vglue 10pt 


{\small\leftskip 25pt\rightskip 25pt{\bf Abstract}\stdspace\theabstract

{\bf AMS Classification}\stdspace\theprimaryclass
\ifx\thesecondaryclass\relax\else; \thesecondaryclass\fi\par
{\bf Keywords}\stdspace \thekeywords\par}
\vglue 7pt
}    
\let\maketitle\makeshorttitle        
%
%

\def\volumenumber#1{\def\thevolumenumber{#1}}
\def\volumeyear#1{\def\thevolumeyear{#1}}
\def\pagenumbers#1#2{\def\startpage{#1}\def\finishpage{#2}}
\def\published#1{\def\publishdate{#1}}
\def\received#1{\def\receiveddate{#1}}
\def\revised#1{\def\reviseddate{#1}}
\let\reviseddate\relax
\volumenumber{X}
\volumeyear{20XX}
\pagenumbers{1}{XXX}
\published{XX Xxxember 20XX}

\long\def\makeagttitle{   
\agt\hfill      
\hbox to 60truept{\vbox to 0pt{\vglue -14truept{\bf [Logo here]}\vss}\hss}
\break
{\small Volume \thevolumenumber\ (\thevolumeyear)
\startpage--\finishpage\nl
Published: \publishdate}

\vglue .2truein

{\parskip=0pt\leftskip 0pt plus 1fil\def\\{\par\smallskip}{\large
\bf\thetitle}\par\medskip}   
\vglue 0.05truein 

%
{\parskip=0pt\leftskip 0pt plus 1fil\def\\{\par}{\sc\theauthors}
\par\medskip}%
 
\vglue 0.03truein 


{\small\leftskip 25truept\rightskip 25truept{\bf Abstract}\stdspace\theabstract

{\bf AMS Classification}\stdspace\theprimaryclass
\ifx\thesecondaryclass\relax\else; \thesecondaryclass\fi\par
{\bf Keywords}\stdspace \thekeywords\par}\vglue 7truept

}   


\def\Addresses{\bigskip
{\small \parskip 0pt \leftskip 0pt \rightskip 0pt plus 1fil \def\\{\par}
\sl\theaddress\par\medskip \rm Email:\stdspace\tt\theemail\par
\ifx\theurl\relax\else\smallskip \rm URL:\stdspace\tt\theurl\par\fi}}

\def\agtart{
\hoffset 14truemm
\voffset 31truemm
\font\phead=cmsl9 scaled 950
\font\pnum=cmbx10 scaled 913
\font\pfoot=cmsl9 scaled 950
\headline{\vbox to 0pt{\vskip -4.5mm\line{\small\phead\ifnum
\count0=\startpage ISSN numbers are printed here
\hfill {\pnum\folio}\else\ifodd\count0\def\\{ }%
\ifx\theshorttitle\relax\thetitle\else\theshorttitle\fi\hfill{\pnum\folio}
\else\def\\{ and }{\pnum\folio}\hfill\ifx\theshortauthors\relax\theauthors
\else\theshortauthors\fi\fi\fi}\vss}}
\footline{\vbox to 0pt{\vglue 0mm\line{\small\pfoot\ifnum\count0=\startpage
Copyright declaration is printed here\hfill\else
\agt, Volume \thevolumenumber\ (\thevolumeyear)\hfill\fi}\vss}}
\let\maketitle\makeagttitle\let\makeshorttitle\makeagttitle}


\def\ifplaintex{\expandafter\ifx\csname documentclass\endcsname\relax}

\def\gtp{{\mathsurround=0pt\it $\cal G\mskip-2mu$eometry \&\ 
$\cal T\!\!$opology $\cal P\!$ublications}}  

\def\recd{{\small Received:\qua\receiveddate\ifx\reviseddate\relax
\else\qquad Revised:\qua\reviseddate\fi\par}} 


\def\lognumber#1{\def\thelognumber{#1}}
\def\volumenumber#1{\def\thevolumenumber{#1}}
\def\volumeyear#1{\def\thevolumeyear{#1}}
\def\papernumber#1{\def\thepapernumber{#1}}
\def\pagenumbers#1#2{\def\startpage{#1}\def\finishpage{#2}}
\def\published#1{\def\publishdate{#1}}

\def\received#1{\def\receiveddate{#1}}
\def\revised#1{\def\reviseddate{#1}}
\def\accepted#1{\def\accepteddate{#1}}
\def\asciititle#1{\def\theasciititle{#1}}

\long\def\asciiabstract#1{\long\def\theasciiabstract{#1}}


\let\\\par\let\thelognumber\relax\let\thevolumenumber\relax
\let\thepapernumber\relax\let\thevolumeyear\relax\let\startpage\relax
\let\finishpage\relax\let\publishdate\relax\let\receiveddate\relax
\let\reviseddate\relax\let\accepteddate\relax\let\theasciititle\relax
\let\theasciiauthors\relax
\let\theasciiabstract\relax

\let\theasciiemail\relax


\ifplaintex
\font\logobig=cmssbx10 scaled 3836
\font\logomed=cmssbx10 scaled 2557
\else
\font\logobig=cmssbx10 scaled 4200
\font\logomed=cmssbx10 scaled 2800
\fi

\long\def\makeagttitle{   
\count0=\startpage
\agt\hfill      
\hbox to 45truept{\vbox to 0pt{\vglue -13truept{\logomed A\kern -.37em{\logobig 
T}\kern -.38em G}\vss}\hss}
\break
{\small Volume \thevolumenumber\ (\thevolumeyear)
\startpage--\finishpage\nl
Published: \publishdate}

\vglue .25truein

{\parskip=0pt\leftskip 0pt plus
1fil\def\\{\par\smallskip}{\Large\bf\thetitle}\par\medskip} \vglue
0.05truein

%
{\parskip=0pt\leftskip 0pt plus 1fil\def\\{\par}{\sc\theauthors}
\par\medskip}%
 
\vglue 0.03truein 


{\small\leftskip 25truept\rightskip 25truept{\bf Abstract}\stdspace\theabstract

{\bf AMS Classification}\stdspace\theprimaryclass
\ifx\thesecondaryclass\relax\else; \thesecondaryclass\fi\par
{\bf Keywords}\stdspace \thekeywords\par}\vglue 7truept

}   

\ifplaintex
\hoffset 14truemm
\voffset 31truemm
\font\phead=cmsl9 scaled 950
\font\pnum=cmbx10 scaled 913
\font\pfoot=cmsl9 scaled 950
\headline{\vbox to 0pt{\vskip -4.5mm\line{\small\phead\ifnum
\count0=\startpage ISSN 1472-2739 (on-line) 1472-2747 (printed)
\hfill {\pnum\folio}\else\ifodd\count0\def\\{ }%
\ifx\theshorttitle\relax\thetitle\else\theshorttitle\fi\hfill{\pnum\folio}
\else\def\\{ and }{\pnum\folio}\hfill\ifx\theshortauthors\relax\theauthors
\else\theshortauthors\fi\fi\fi}\vss}}
\footline{\vbox to 0pt{\vglue 0mm\line{\small\pfoot\ifnum\count0=\startpage
\copyright\ \gtp\hfill\else
\agt, Volume \thevolumenumber\ (\thevolumeyear)\hfill\fi}\vss}}
\else
\headsep 23pt
\footskip 35pt
\hoffset -4truemm
\voffset 12.5truemm
\font\lhead=cmsl9 scaled 1050
\font\lnum=cmbx10 
\font\lfoot=cmsl9 scaled 1050
\makeatletter
\def\@oddhead{{\small\lhead\ifnum\count0=\startpage ISSN 1472-2739 
(on-line) 1472-2747 (printed)\hfill {\lnum\number\count0}\else\ifodd\count0
\def\\{ }\ifx\theshorttitle\relax \thetitle \else\theshorttitle\fi\hfill
{\lnum\number\count0}\else\def\\{ and }{\lnum\number\count0}
\hfill\ifx\theshortauthors\relax 
\theauthors\else\theshortauthors\fi\fi\fi}}\def\@evenhead{\@oddhead}
\def\@oddfoot{\small\lfoot\ifnum\count0=\startpage\copyright\ \gtp\hfill\else
\agt, Volume \thevolumenumber\ (\thevolumeyear)\hfill\fi}
\def\@evenfoot{\@oddfoot}
\makeatother
\fi
\let\maketitlepage\makeagttitle
\let\makeshorttitle\maketitlepage
\let\maketitle\maketitlepage


\newwrite\gtoutfile
\long\gdef\makeheadfile{  
{\def\\{, }\def\s{ }
\immediate\openout\gtoutfile head.xxx
\immediate\write\gtoutfile{To: math@arxiv.org}
\immediate\write\gtoutfile{Subject: put OR rep NNNNN:ppppp}
\immediate\write\gtoutfile{--text follows this line--}
\immediate\write\gtoutfile{Proxy-for: \ifx\theasciiauthors\relax
\theauthors\else\theasciiauthors\fi\s<\ifx\theasciiemail\relax\theemail\else\theasciiemail\fi>}
\immediate\write\gtoutfile{\noexpand\\}
\immediate\write\gtoutfile{Authors: \ifx\theasciiauthors\relax
\theauthors\else\theasciiauthors\fi}
{\def\\{ }\immediate\write\gtoutfile{Title: \ifx\theasciititle\relax
\thetitle\else\theasciititle\fi}}
\immediate\write\gtoutfile{Subj-class: GT or SG, GR etc}
\immediate\write\gtoutfile{MSC-class: \theprimaryclass\ifx\thesecondaryclass\relax\else, \thesecondaryclass\fi}
\immediate\write\gtoutfile{Journal-ref: Algebr. Geom. Topol. \thevolumenumber\s
(\thevolumeyear) \startpage-\finishpage}
\immediate\write\gtoutfile{Comments: Published by Algebraic and
Geometric Topology at}
\immediate\write\gtoutfile{\s\s\s  http://www.maths.warwick.ac.uk/agt/AGTVol\thevolumenumber/agt-\thevolumenumber-\thepapernumber.abs.html}
\immediate\write\gtoutfile{\noexpand\\}
\immediate\write\gtoutfile{}
\ifx\theasciiabstract\relax
\immediate\write\gtoutfile{\theabstract}\else
\immediate\write\gtoutfile{\theasciiabstract}\fi
\immediate\write\gtoutfile{}
\immediate\write\gtoutfile{\noexpand\\}
\immediate\write\gtoutfile{}
\immediate\closeout\gtoutfile}}  

\def\maketitlepage{\makeagttitle\makeheadfile}
\let\makeshorttitle\maketitlepage
\let\maketitle\maketitlepage

\lognumber{4}
\volumenumber{3}
\volumeyear{2003}
\papernumber{4}
\published{8 February 2003}
\pagenumbers{103}{116}
\received{25 April 2002}
\revised{15 January 2003}
\accepted{7 February 2003}

\input xypic
\input rlepsf
\input newinsert

\def\rtimes{{\msb \char'157 }}
\def\simeq{{\msa \char'167}}

\reflist
\key{AFR} {\bf C.S. Aravinda}, {\bf F.T. Farrell} and {\bf S.K.
Roushon}, {\it Algebraic $K$-theory of pure braid groups},
Asian J. Math. 4 (2000) 337--344 

\key{DM} {\bf J.D. Dixon} and {\bf B. Mortimer} 
{\it Permutation groups}, Grad. Texts in Math. {\bf 163}, Springer, New
York (1996)

\key{E} {\bf D.B.A. Epstein}, {\it Ends}, from:  ``Topology of
3-manifolds and related topics (Proc. The Univ.
of Georgia Institute, 1961)'', Prentice-Hall, Englewood Cliffs, N.J.
(1962) 110--117

\key{FJ}
{\bf F.T. Farrell} and {\bf L.E. Jones}, {\it Isomorphism conjectures in
algebraic $K$-theory}, J. Amer. Math. Soc. 6 (1993) 249--297

\key{FL} {\bf F.T. Farrell} and {\bf P.A. Linnell},
{\it K-Theory of solvable groups}, Proc. London Math. Soc., to appear,
{\tt arXiv:math.KT/0207138}

\key{FR} {\bf F.T. Farrell} and {\bf S.K. Roushon},
{\it The Whitehead groups of braid groups vanish}, Internat. Math. Res.
Notices, no. 10 (2000) 515--526

\key{L} {\bf B. Leeb}, {\it 3-manifolds with(out) metrics of
nonpositive curvature}, Invent. Math. 122 (1995) 277--289

\key{Q} {\bf F. Quinn}, {\it Ends of maps. II}, Invent. Math. 68
(1982) no. 3 353--424

\key{R} {\bf S.K. Roushon}, {\it Fibered isomorphism conjecture for
complex manifolds}, Tata Institute, preprint, {\tt arXiv:math.GT/0209119}

\key{S} {\bf John Stallings}, {\it Topologically unrealizable
automorphisms of free groups}, Proc. Amer. Math. Soc. 84 (1982) 21--24
\endreflist

\title{$K$-theory of virtually poly-surface groups}
\asciititle{K-theory of virtually poly-surface groups}

\author{S.K. Roushon}
\address{School of Mathematics, Tata Institute\\Homi Bhabha Road, Mumbai
400 005, India.}
\email{roushon@math.tifr.res.in}  
\url{http://www.math.tifr.res.in/\char'176roushon/paper.html}

\abstract In this paper we generalize the notion of strongly poly-free
group to a larger class of groups, we call them {\it strongly
poly-surface} groups and prove that the Fibered Isomorphism Conjecture of
Farrell and Jones corresponding to the stable topological pseudoisotopy
functor is true for any virtually strongly poly-surface group.
A consequence is that the Whitehead group of a torsion free subgroup of
any virtually strongly poly-surface group vanishes. \endabstract

\asciiabstract{In this paper we generalize the notion of strongly poly-free
group to a larger class of groups, we call them strongly
poly-surface groups and prove that the Fibered Isomorphism Conjecture of
Farrell and Jones corresponding to the stable topological pseudoisotopy
functor is true for any virtually strongly poly-surface group.
A consequence is that the Whitehead group of a torsion free subgroup of
any virtually strongly poly-surface group vanishes.}

\primaryclass{19B28, 19A31, 20F99, 19D35}
\secondaryclass{19J10}

\keywords{Strongly poly-free groups, poly-closed surface groups, Whitehead
group, fibered isomorphism conjecture}
\maketitle

\section{Introduction}

We generalize the class of strongly poly-free groups which was introduced  
in [\ref{AFR}]. 

\rk{Definition 1.1} A discrete group $\Gamma$ is called 
{\it strongly poly-surface} if there exists a finite filtration of
$\Gamma$ by subgroups: $1=\Gamma_0\subset \Gamma_1\subset \cdots \subset
\Gamma_n=\Gamma$ such that the following conditions are satisfied:

\items
\item{(1)} $\Gamma_i$ is normal in $\Gamma$ for each $i$.

\item{(2)} $\Gamma_{i+1}/\Gamma_i$ is isomorphic to the fundamental group
of a surface.

\item{(3)} for each $\gamma\in \Gamma$ and $i$ there is a 
surface $F$ such that $\pi_1(F)$ is isomorphic to
$\Gamma_{i+1}/\Gamma_i$ and either (a) $\pi_1(F)$ is finitely generated or
(b) $\pi_1(F)$ is infinitely generated and $F$ has one end. Also there is
a diffeomorphism $f:F\to F$ such that the induced outer 
automorphism $f_{\#}$ of $\pi_1(F)$ is equal to $c_\gamma$ in 
$Out(\pi_1(F))$, where $c_\gamma$ is the outer automorphism of
$\Gamma_{i+1}/\Gamma_i\ \simeq \ \pi_1(F)$ induced by the conjugation
action on $\Gamma$ by $\gamma$.
\enditems

In such a situation we say that the group $\Gamma$ has {\it rank} $\leq
n$. 
\endrk

Note that in the definition of strongly poly-free
group we demanded that the groups $\Gamma_{i+1}/\Gamma_i$ be finitely
generated free groups. On the other hand in the definition of strongly
poly-surface group, $\Gamma_{i+1}/\Gamma_i$ can be the fundamental
group of any surface other than the surfaces with infinitely generated
fundamental groups and with more than one topological ends. We even allow
a class of surfaces with infinitely generated fundamental group. Also we
remark that if the groups in $(2)$ are fundamental groups of closed
surfaces then the condition $(3)$ is always satisfied. This follows from
the well-known fact that any automorphism of the fundamental group of a
closed surface is induced by a
diffeomorphism of the surface. However this fact is very rarely true
for surfaces with nonempty boundary ([\ref{S}]). Thus the class of
strongly
poly-surface groups contains a class of poly-closed surface groups. Here
recall that given a class of groups $\cal G$, a group $\Gamma$ is called
poly-$\cal G$ if $\Gamma$ has a filtration by subgroups $1=\Gamma_0\subset
\Gamma_1\subset \cdots \subset \Gamma_n=\Gamma$ such that $\Gamma_i$
is normal in $\Gamma_{i+1}$ and $\Gamma_{i+1}/\Gamma_i\in {\cal G}$ for
each $i$. And a group is called virtually poly-$\cal G$ if it has a
normal subgroup $G\in {\cal G}$ of finite index. For a group $G$, 
by `poly-$G$' we will mean `poly-$\cal G$', where $\cal G$ consists of
$G$ only.

In [\ref{AFR}] we proved that the Whitehead group of any strongly
poly-free
group vanishes. Generalizing this result the Fibered Isomorphism
Conjecture (FIC) corresponding to the stable topological pseudoisotopy
functor ([\ref{FJ}]) was proved for any virtually strongly poly-free group
in [\ref{FR}]. In this paper we prove FIC for any virtually
strongly poly-surface group. The Main Lemma in the next section is the
crucial result which makes this generalization possible. The key idea to
prove the Main Lemma is that, except for three closed surfaces, the
covering space corresponding to the commutator subgroup of the
fundamental group of all other closed surfaces have one topological end.  

Below we recall the Fibered Isomorphism Conjecture in brief. For details
about this conjecture see [\ref{FJ}]. Here we follow the formulation given
in [\ref{FL}, appendix]. 

Let $\cal S$ denote one of the three functors from the category of
topological spaces to the category of spectra: (a) the stable topological
pseudo-isotopy functor ${\cal P}()$; (b) the algebraic $K$-theory functor 
${\cal K}()$; (c) and the $L$-theory functor ${\cal L}^{-\infty}()$. 

Let $\cal M$ be the category of continuous surjective maps. The
objects of $\cal M$ are continuous surjective maps $p:E\to B$ between
topological spaces $E$ and $B$. And a morphism between two maps $p:E_1\to
B_1$ and $q:E_2\to B_2$ is a pair of continuous maps $f:E_1\to E_2$,
$g:B_1\to B_2$ such that the following diagram commutes. 
$$\diagram
E_1 \rto^f \dto^p & E_2 \dto^q\\
B_1 \rto^g &B_2\enddiagram$$
There is a functor defined by Quinn [\ref{Q}] from $\cal M$ to the category
of $\Omega$-spectra which associates to the map $p$ the spectrum ${\Bbb 
H}(B, {\cal S}(p))$ with the property that ${\Bbb H}(B, {\cal S}(p))={\cal
S}(E)$ when $B$ is a single point. For an explanation of ${\Bbb H}(B,
{\cal S}(p))$ see [\ref{FJ}, section 1.4]. Also the map ${\Bbb H}(B,
{\cal S}(p))\to {\cal S}(E)$ induced by the morphism: id$:E\to E$; 
$B\to *$ in the category $\cal M$ is called the Quinn assembly map.

Let $\Gamma$ be a discrete group and $\cal E$ be a $\Gamma$ space which
is universal for the class of all virtually cyclic subgroups of $\Gamma$ 
and denote ${\cal E}/\Gamma$ by $\cal B$. For definition of universal 
space see [\ref{FJ}, appendix]. Let $X$ be a space on which $\Gamma$ acts
freely and properly discontinuously and $p:X\times_{\Gamma} {\cal E}\to
{\cal E}/{\Gamma}={\cal B}$ be the map induced by the projection onto the 
second factor of $X\times {\cal E}$. 

The Fibered Isomorphism Conjecture states that the map $${\Bbb
H}({\cal B}, {\cal S}(p))\to {\cal S}(X\times_{\Gamma} {\cal E})={\cal
S}(X/\Gamma)$$ is an (weak) equivalence of spectra. The equality in the
above display is induced by the map $X\times_{\Gamma}{\cal E}\to X/\Gamma$
and using the fact that $\cal S$ is homotopy invariant.  

Let $Y$ be a connected $CW$-complex and $\Gamma \ \simeq \ \pi_1(Y)$. Let
$X$ be the universal cover $\tilde Y$ of $Y$ and the action of $\Gamma$ on
$X$ is the action by group of covering transformation. If we take an
aspherical $CW$-complex $Y'$ with $\Gamma \ \simeq \ \pi_1(Y')$ and $X$ is
the universal cover $\tilde Y'$ of $Y'$ then by [\ref{FJ}, corollary
2.2.1] if the FIC is true for the space $\tilde Y'$ then it is true for
$\tilde Y$ also. Thus whenever we say that FIC is true for a discrete
group $\Gamma$ or for the fundamental group $\pi_1(X)$ of a space $X$ we
shall mean it is true for the Eilenberg-MacLane space $K(\Gamma, 1)$ or
$K(\pi_1(X), 1)$ and for the functor ${\cal S}()$. 

Throughout this paper we consider only the stable topological
pseudo-isotopy functor; that is the case when ${\cal S}()={\cal P}()$. 
And by FIC we mean FIC for ${\cal P}()$.

The main theorem of this article is the following.

\proclaim{Main Theorem} Let $\Gamma$ be a virtually strongly poly-surface
group. Then the Fibered Isomorphism Conjecture is true for
$\Gamma$.\endproc

Recall that if FIC is true for a torsion free group $G$ then
$Wh(G)={\tilde K}_0({\Bbb Z}G)=K_{-i}({\Bbb Z}G)=0$ for all $i\geq 1$. 
A proof of this fact is given in several places, e.g., see [\ref{FR}] 
or [\ref{FL}].

Hence we have the following corollary.

\proclaim{Corollary 1.2} Let $G$ be a torsion free subgroup of a virtually
strongly poly-surface group. Then $Wh(G)={\tilde K}_0({\Bbb
Z}G)=K_{-i}({\Bbb Z}G)=0$ for all $i\geq 1$.\endproc

\section{Proof of the Main Theorem}

The proof of the Main Theorem appears at the end of this section. Before
that we state some known results about the Fibered Isomorphism Conjecture
and prove the Main Lemma and some propositions. Apart from being 
crucial ingredients to the proof of the Main Theorem the Main Lemma and the
propositions are also of independent interest. 

Recall that the FIC is true for any finite group and for abelian groups
([\ref{FL}, lemma 2.7]).

\proclaim{Lemma A} {\rm([\ref{FJ}, theorem A.8])}\qua If the FIC is true for a 
discrete group $\Gamma$ then it is true for any subgroup of
$\Gamma$.\endproc

Before we state the next lemma let us recall the following group theoretic
definition. Let $G$ and $H$ be two groups. Assume $G$ is finite. Then
$H\wr G$ denotes the wreath
product with respect to the regular action of $G$ on $G$. Recall that
actually $H\wr G \ \simeq \ H^G\rtimes G$ where $H^G$ is product of $|G|$
copies of $H$ indexed by $G$ and $G$ acts on the product via the regular
action of $G$ on $G$. An easily checked fact is that if $G_1$ is another
finite group then $H^{G_1}\wr G$ is a subgroup of $H\wr (G_1\times G)$.
This fact will be used throughout the paper. Another
fact we will be using is that for any two groups $A$ and $B$ the group
$(A\times B)\wr G$ is a subgroup of $(A\wr G)\times (B\wr G)$. 

The Algebraic Lemma from [\ref{FR}] says the following.

\proclaim{Algebraic Lemma} If $G$ is an extension of a group
$H$ by a finite group $K$ then $G$ is a subgroup of $H\wr K$.\endproc
 
This lemma is also proved in [\ref{DM}, theorem 2.6A].

\proclaim{Lemma B} {\rm[\ref{R}]}\qua Let $\Gamma$ be an extension of the
fundamental group $\pi_1(M)$ of a closed nonpositively curved 
Riemannian manifold or a compact surface (may be with nonempty 
boundary) $M$ by a finite group $G$ then FIC is true for
$\Gamma$. Moreover FIC is true for the wreath product $\Gamma\wr G$. 
\endproc

\prf Let us consider the closed case first. By the Algebraic Lemma
we have an embedding of $\Gamma$ in the wreath product $\pi_1(M)\wr
G$. Let $U=M\times \cdots \times M$ be the $|G|$-fold product of $M$.
Then $U$ is a closed nonpositively curved Riemannian manifold. By
[\ref{FR}, fact 3.1] it follows that FIC is true for
$\pi_1(U)\rtimes G \ \simeq \ (\pi_1(M))^G\rtimes G \ \simeq \ \pi_1(M)\wr
G$. Lemma A now proves that FIC is true for $\Gamma$.

If $M$ is a compact surface with nonempty boundary then $\pi_1(M) < 
\pi_1(N)$ where $N$ is a closed nonpositively curved surface. Hence
$\Gamma < \pi_1(M)\wr G < \pi_1(N)\wr G$. Using Lemma A and the previous 
case we complete the proof.\endprf

The above Lemma is also true if $M$ is a compact irreducible $3$-manifold
with nonempty incompressible boundary and the boundary components are
torus or Klein bottle. Indeed in this situation by theorem 3.2 and 3.3
from [\ref{L}] the interior of $M$ supports a complete nonpositively
curved Riemannian metric so that near the boundary the metric is a product
metric. Hence the double of $M$ will support a nonpositively curved
metric and we argue as in the case of compact surface to deduce the
following Corollary.

\proclaim{Corollary B} Let $M_1,\cdots, M_k$ be compact irreducible
$3$-manifolds with incompressible boundary which has either
torus or Klein bottle as components. Then FIC is true for
$(\pi_1(M_1)\times\cdots \times \pi_1(M_n))\wr G$ for any finite
group $G$.\endproc

\proclaim{Lemma C} {\rm([\ref{FJ}, proposition 2.2])}\qua Let $f:G\to H$ be a 
surjective homomorphism. Assume that the FIC is true for $H$ and for
$f^{-1}(C)$ for all virtually cyclic subgroup $C$ of $H$ (including
$C=1$). Then FIC is true for $G$.\endproc

We will use Lemma A, Lemma C and the Algebraic Lemma 
throughout the paper, sometimes even without referring to them. 

We now recall a well-known fact from $2$-dimensional real manifold theory.

\proclaim{Lemma D} Let $\Gamma$ be a finitely generated nonabelian free
group. Then $\Gamma$ is isomorphic to the fundamental group of a compact
surface (with nonempty boundary).\endproc

\proclaim{Lemma E} Let $\Gamma$ be the fundamental group of a surface then
FIC is true for $\Gamma\wr G$ for any finite group $G$.\endproc

\prf If $\Gamma$ is finitely generated then $\Gamma$ is the
fundamental group of a compact surface and hence  
the lemma follows from Lemma B. In the infinitely generated
case $\Gamma\wr G \ \simeq \ \lim_{i\to \infty} (\Gamma_i\wr G)$ where
each $\Gamma_i$ is a finitely generated nonabelian free group. By Lemma B,
Lemma D and Theorem F (see below) the proof is complete.\endprf

We quote the following theorem of Farrell and Linnell which will be
used throughout the paper.

\proclaim{Theorem F}{\rm([\ref{FL}, theorem 7.1])}\qua 
Let $I$ be a directed set, and let $\Gamma_n$, $n\in I$ be a directed
system of groups with $\Gamma=\lim_{n\in I}$; i.e., $\Gamma$ is the
direct limit of the groups $\Gamma_n$. If each group $\Gamma_n$
satisfies FIC, then $\Gamma$ also satisfies FIC.\endproc

We will also use proposition 2.4 from [\ref{FJ}] frequently, sometime
without referring to it. This result says that FIC is true for any
virtually poly-$\Bbb Z$ group. 

We need the following crucial proposition to prove the Main Theorem.

\proclaim{Proposition 2.1} Let $S$ be a surface. If $\pi_1(S)$ is
infinitely generated then assume $S$ has one topological end. Let $f$ be a
diffeomorphism of $S$. Then the group $\pi_1(S)\rtimes {\Bbb Z}$ satisfies
the FIC. Here, up to conjugation, the action of a generator of ${\Bbb
Z}$ on the group $\pi_1(S)$ is induced by the diffeomorphism
$f$.\endproc 

\prf There are two cases according as $S$ is compact or not. 

If $S$ is compact with nonempty boundary then
$\pi_1(S)\rtimes {\Bbb Z}$ is the fundamental group of a compact
irreducible $3$-manifold $M$ with torus or Klein bottle as boundary
component. If $\pi_1(S)=1$ then there is nothing to prove, otherwise the
boundary components of $M$ will be incompressible. Hence Corollary B
proves this case.

So assume that either $S$ is closed or a noncompact surface. Note that if
the fundamental group is finitely generated free then by Lemma D it falls
in the previous case.

Let us consider the closed case first.
This case is contained in the following Lemma which was proved in
[\ref{R}] in the case when the fiber is orientable. Here we give a proof 
for the general situation.

\proclaim{Main Lemma} Let $M^3$ be a closed $3$-dimensional 
manifold which is the total space of a fiber bundle projection $M^3\to
{\Bbb S}^1$ with fiber $F$ such that $b_1(F)\neq 1$. Then FIC is true for
$\pi_1(M)$.\endproc  

\prf 
The following exact sequence is obtained from the long exact homotopy
sequence of the fibration $M\to {\Bbb S}^1$. $$1\to \pi_1(F)\to
\pi_1(M)\to \pi_1({\Bbb S}^1)\to 1$$
Let $[A,A]$ denotes the commutator subgroup of the group
$A$. Then we have $$1\to [\pi_1(F), \pi_1(F)]\to \pi_1(F)\to H_1(F, {\Bbb
Z})\to 1$$ Let $t$ be a generator of $\pi_1({\Bbb S}^1)$. Since
$[\pi_1(F), \pi_1(F)]$ is a characteristic subgroup of $\pi_1(F)$ the
action (induced by the monodromy) of $t$ on $\pi_1(F)$ leaves $[\pi_1(F),
\pi_1(F)]$ invariant. Thus we have another exact sequence 
$$1\to [\pi_1(F), \pi_1(F)]\to \pi_1(F)\rtimes
\langle t\rangle \to H_1(F,
{\Bbb Z})\rtimes \langle t\rangle \to 1$$ Which reduces to the sequence
$$1\to [\pi_1(F), \pi_1(F)]\to \pi_1(M)\to H_1(F, {\Bbb
Z})\rtimes \langle t\rangle \to 1$$
We would like to apply Lemma C to this exact sequence. 

Now we have two cases according as the fiber is orientable 
or nonorientable. Let us first consider the orientable fiber case. If the
fiber is ${\Bbb S}^2$ or ${\Bbb T}^2$ then $\pi_1(M)$ is poly-$\Bbb Z$ and
hence FIC is true for $\pi_1(M)$. So assume that the fiber has
genus $\geq 2$.

Clearly the group $H_1(F, {\Bbb Z})\rtimes \langle t\rangle $ is
poly-$\Bbb Z$. Hence FIC is true for $H_1(F, {\Bbb
Z})\rtimes \langle t\rangle$. Let $C$ be a virtually cyclic subgroup of 
$H_1(F, {\Bbb Z})\rtimes \langle t\rangle$. Let $p:\pi_1(M)\to H_1(F,
{\Bbb Z})\rtimes \langle t\rangle$ be the above
surjective homomorphism. We will show that 
the FIC is true for $p^{-1}(C)$. Note that $C$ is either trivial or
infinite cyclic. 

{\bf Case $C=1$}\quad In this case we have that $p^{-1}(C)$ is 
a nonabelian free group and hence is the fundamental group of a surface.
Lemma E proves this case. 

{\bf Case $C\neq 1$}\quad We have $p^{-1}(C) \ \simeq \ [\pi_1(F),
\pi_1(F)]\rtimes \langle s\rangle $ where $s$ is a generator of $C$. Let
$\tilde F$ be the covering space of $F$
corresponding to the commutator subgroup $[\pi_1(F), \pi_1(F)]$. As $F$
has first Betti number $\geq 2$ the group $H_1(F, {\Bbb Z})$ has only one
end. Also $H_1(F, {\Bbb Z})$ is the group of covering transformations of
the regular covering $\tilde F\to F$. Since $F$ is compact, the manifold
$\tilde F$ has one topological end (see [\ref{E}]). Figure 1 
describes $\tilde F$.

We write the manifold $\tilde F$ as the union of compact submanifolds. 
As $\tilde F$ has one end there is a connected compact submanifold $M_0$
of $\tilde F$ so that the complement $\tilde F-M_0$ has one
connected component and
for any other connected compact submanifold $M$ containing $M_0$ the
complement $\tilde F-M$ also has one component. Consider a sequence $M_i$
of compact submanifolds of $\tilde F$ with the following properties.

\items
\item{(1)} each $M_i$ has one boundary component

\item{(2)} each $M_i$ has the same property as  $M_0$

\item{(3)} $\tilde F=\cup_iM_i$\quad and 

\item{(4)} $M_0\subset M_1\subset\cdots$ .
\enditems

\figure
\relabelbox\small
\cl{\epsfxsize 3in \epsfbox{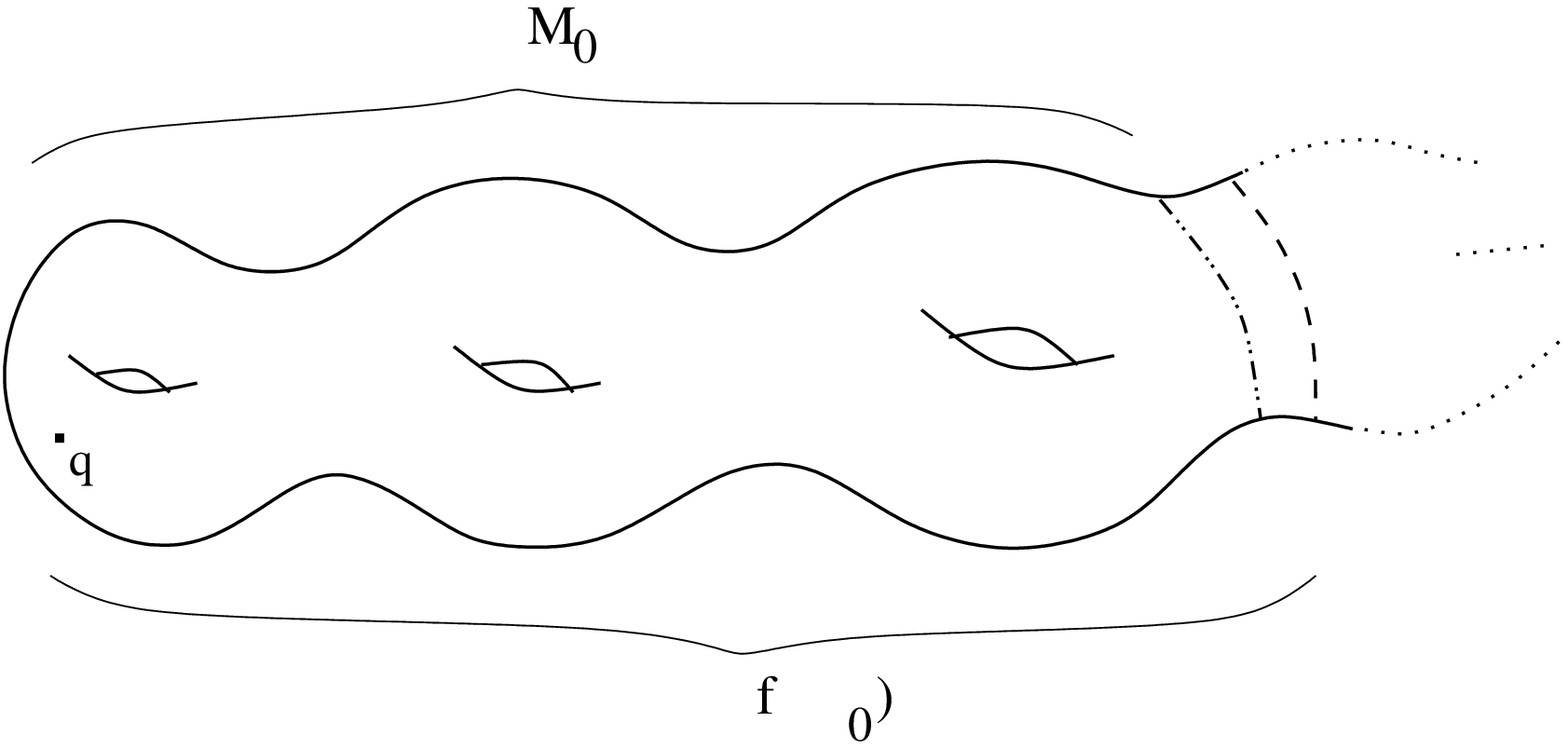}}
\relabel {q}{$q$}
\relabel {M}{$M_0$}
\relabel {f}{$f(M_0)$}
\endrelabelbox
\endfigure

Note that the monodromy diffeomorphism of $F$ lifts to a
quadiffeomorphism of $\tilde F$ which in turn, up to conjugation, induces the
action of $t$ on $[\pi_1(F), \pi_1(F)]$ and also $H_1(F, {\Bbb Z})$ is
the group of covering transformation of $\tilde F\to F$. Also, the
induced action of $t$ on $H_1(F, {\Bbb Z})$ is given by $t(s)=\tilde
f\circ s\circ {\tilde f}^{-1}$, where $\tilde f:\tilde F\to \tilde F$ is a
lift of the monodromy diffeomorphism and $s\in H_1(F, {\Bbb Z})$ acts on
$\tilde F$ as a covering transformation. From this observation it follows
that, up to conjugation, the action of $s$ on $[\pi_1(F), \pi_1(F)]$ is
induced by a diffeomorphism (say $f$) of $\tilde F$. Indeed, if $s=(s_1,
t^k)\in H_1(F, {\Bbb Z})\rtimes \langle t\rangle$ then $f=s_1\circ {\tilde
f}^k:\tilde F\to \tilde F$.

Note that $f$ is also a lift of a diffeomorphism (say $f_1$) of $F$. If
$f_1$ is isotopic to a pseudo-Anosov diffeomorphism then by Thurston's
theorem $p^{-1}(C)$ is a subgroup of the fundamental group of a closed
hyperbolic $3$-manifold, namely the mapping torus of $f_1$. Hence FIC is
true for $p^{-1}(C)$. So we can assume that $f_1$ is isotopic to either a
finite order diffeomorphism or to a reducible one. Hence there exists in
$\tilde F$ simple closed curves so that cutting along them produces a
filtration of $\tilde F$ with properties $(1)$ to $(4)$.

Using properties $(1)$ to $(4)$ and that $f_1$ is either of finite
order or reducible, it is now easy to see that each $f(M_i)$ is obtained
from $M_i$ by attaching an annulus to the boundary component of $M_i$. So,
we can isotope $f$ so that $f(M_i)=M_i$ for each $i$. Thus we have a
filtration $\pi_1(M_0, q)<\pi_1(M_1, q)<\cdots$ of $[\pi_1(F), \pi_1(F)]$
by finitely generated  free subgroups so that the action of $s$ on
$[\pi_1(F), \pi_1(F)]$ respects this filtration and each
$\pi_1(M_i)\rtimes \langle s\rangle$ is the fundamental group of a Haken
$3$-manifold $N^s_i$ with nonempty incompressible boundary. Indeed,
$N^s_i$ is diffeomorphic to the mapping torus of the restriction of $f$ to
$M_i$. 

Hence FIC is true for $\pi_1(N^s_i)$ by Corollary B. From above we
also get that $[\pi_1(F), \pi_1(F)]\rtimes \langle s\rangle \ \simeq \ 
\lim_{i \to \infty } \pi_1(N^s_i)$. Using Theorem F we conclude that FIC
is true for $[\pi_1(F), \pi_1(F)]\rtimes \langle s\rangle $.

This completes the proof of the Main Lemma in the orientable fiber case.

From the above proof we get the following Lemma which is true for
nonorientable $F$ also.

\proclaim{Lemma 2.2} Let $F$ be a closed surface of genus $\geq 2$ and
$\tilde F$ be the covering of $F$ corresponding to the commutator subgroup
of $\pi_1(F)$. Let $f$ be a diffeomorphism of $F$. Then $[\pi_1(F),
\pi_1(F)]\rtimes \langle s\rangle \ \simeq \ \lim_{i \to \infty
}\pi_1(N^s_i)$ where $N^s_i$ are compact Haken $3$-manifolds with nonempty
incompressible boundary and up to
conjugation the action of $s$ on $[\pi_1(F), \pi_1(F)]$ is induced by the
lift of $f$ to $\tilde F$. Moreover each $\pi_1(N^s_i)$ is a
subgroup of the fundamental group of a closed nonpositively curved
Riemannian manifold $M^s_i$.\endproc

Next we deal with the case when the fiber is nonorientable. In this
situation $H_1(F, {\Bbb Z})$ has torsion element. Nevertheless
$H_1(F, {\Bbb Z})\rtimes\langle t\rangle$ is a virtually poly-$\Bbb Z$ 
group and hence FIC is true for this group. Thus we can apply Lemma C to the exact sequence.  
$$1\to [\pi_1(F), \pi_1(F)]\to \pi_1(M)\to H_1(F, {\Bbb
Z})\rtimes \langle t\rangle \to 1$$
If $F$ is the projective plane then $\pi_1(M)$ is virtually
infinite cyclic and FIC is true for this group. Since the fiber is not
the Klein bottle we assume that genus of $F$ is $\geq 2$.

Again we have two cases.

\noindent
{\bf $C$ is finite}\quad We have 
$$p^{-1}(C) < ([\pi_1(F), \pi_1(F)])\wr C$$
Hence FIC is true for $p^{-1}(C)$ by Lemma E.

\noindent
{\bf $C$ is infinite}\quad 
Let $C_1$ be an infinite cyclic subgroup of $C$ of finite index. As $C_1$
is of finite index we can assume that $C_1$ is normal in $C$. We have the
following exact sequences. $$1\to p^{-1}(C_1)\to p^{-1}(C)\to G\to 1$$ and
$$1\to [\pi_1(F), \pi_1(F)]\to p^{-1}(C_1)\to C_1\to 1$$
Here $G$ is a finite group. Let $C_1$ be generated by $s$. Then we get 
$$p^{-1}(C) < ([\pi_1(F), \pi_1(F)]\rtimes \langle s\rangle )\wr G$$
As in the orientable case, up to conjugation, the action of $s$
on $[\pi_1(F), \pi_1(F)]$ is induced by a diffeomorphism of $\tilde F$.
Also recall that genus of $F$ is $\geq 2$. Hence Lemma 2.2 is applicable.
Thus we get $$([\pi_1(F), \pi_1(F)]\rtimes \langle s\rangle )\wr G \ 
\simeq \ \lim_{i \to \infty }(\pi_1(N^s_i)\wr G)$$
Lemma B together with Theorem F complete the proof in this
case.\endprf

To complete the proof of Proposition 2.1 we need to consider the case 
when $\pi_1(S)$ is infinitely generated and $S$ has one end. We use Lemma
2.2 to deduce that $\pi_1(S)\rtimes {\Bbb Z} \ \simeq \ \pi_1(S)\rtimes
\langle
t\rangle \ \simeq \ \lim_{i \to \infty }\pi_1(N^t_i)$. Now
apply Corollary B and Theorem F to complete the proof of the
proposition.\endprf

The proposition below is an application of the method of the proof of the
Main Lemma.

\proclaim{Proposition 2.3} Let $M$ be as in the Main Lemma. Then FIC is
true for $\pi_1(M)\wr G$ for any finite group $G$.\endproc

\prf Recall the following exact sequence. $$1\to [\pi_1(F),
\pi_1(F)]\to \pi_1(M)\to H_1(F, {\Bbb Z})\rtimes \langle t\rangle \to 1$$
If $F$ is the $2$-sphere or the projective plane then $\pi_1(M)\wr G$ is
virtually abelian and hence FIC is true by [\ref{FL}, lemma 2.7]. So
assume $F$ is not the $2$-sphere or the Klein bottle or the projective
plane.

Taking wreath product with $G$ the above exact sequence gives the
following. $$1\to ([\pi_1(F), \pi_1(F)])^G\to \pi_1(M)\wr G\to (H_1(F,
{\Bbb Z})\rtimes \langle t\rangle)\wr G\to 1$$
Recall that $(H_1(F, {\Bbb Z})\rtimes \langle t\rangle)\wr G$ is virtually
poly-$\Bbb Z$ and hence FIC is true for $(H_1(F, {\Bbb Z})\rtimes \langle
t\rangle)\wr G$. Applying Lemma C twice and noting that FIC is true for
free abelian groups, it is easy to show that if the FIC is true for two
torsion free group then it is true for the product of the two groups also. 
Thus by Theorem F and Lemma E it follows that FIC is true for $([\pi_1(F),
\pi_1(F)])^G$. 

Let $Z$ be a virtually cyclic subgroup of $(H_1(F, {\Bbb Z})\rtimes
\langle t\rangle)\wr G$. If $Z$ is finite then $$p^{-1}(Z) < 
([\pi_1(F), \pi_1(F)])^G\wr Z < ([\pi_1(F),
\pi_1(F)])\wr (G\times Z)$$
Here $p$ is the surjective homomorphism $\pi_1(M)\wr G\to (H_1(F, {\Bbb
Z})\rtimes \langle t\rangle)\wr G$.

Now Lemma E applies on the right hand side group to show that FIC is true
for $p^{-1}(Z)$. If $Z$ is infinite then let $Z_1$ be the intersection of
$Z$ with the torsion free part of $(H_1(F, {\Bbb Z})\rtimes \langle
t\rangle)^G$. Hence $Z_1 \ \simeq \ \langle u\rangle$ is infinite cyclic
normal subgroup of $Z$ of finite index. Once again we appeal to the
Algebraic Lemma to get $$p^{-1}(Z) < (p^{-1}(Z_1))\wr {Z/Z_1} \ \simeq \ 
(([\pi_1(F),
\pi_1(F)])^G\rtimes \langle u\rangle )\wr Z/Z_1$$ $$ \ \simeq \  
(([\pi_1(F),
\pi_1(F)]\times [\pi_1(F), \pi_1(F)]\cdots\times [\pi_1(F), 
\pi_1(F)])\rtimes \langle u\rangle )\wr
Z/Z_1=H\text{(say)}$$
In the above display there are $|G|$ number of factors of $[\pi_1(F),
\pi_1(F)]$. Note that the action of $u$ on $([\pi_1(F), \pi_1(F)])^G$ 
is factorwise. Let us denote the restriction of the action of $u$ on the
$j$-th factor of $([\pi_1(F), \pi_1(F)])^G$ by $u_j$. By Lemma 2.2 we get
$$H < (\lim_{i \to \infty }(\pi_1(M^{u_1}_i)\times
\pi_1(M^{u_2}_i)\times\cdots\times\pi_1(M^{u_{|G|}}_i))\wr Z/Z_1 \ \simeq
\ \lim_{i \to \infty }(\pi_1(M_i)\wr Z/Z_1)$$ where $M^{u_j}_i$ and hence 
$M_i=M^{u_1}_i\times M^{u_2}_i\times\cdots\times M^{u_{|G|}}_i$ are closed
nonpositively curved Riemannian manifolds. Using Lemma B and Theorem F we
complete the proof of the Proposition.\endprf

The following corollary is a consequence of Proposition 2.3.

\proclaim{Corollary 2.4} Let $M_i$ for $i=1,\cdots , k$ be $3$-manifolds  
with the same property as $M$ in the Main Lemma. Then FIC is true for
$(\pi_1(M_1)\times\cdots \times\pi_1(M_k))\wr G$ for any finite group
$G$.\endproc 

\prf For the proof of the Corollary just note that if $A$ and $B$
be two groups and $G$ is another group acting regularly on itself then
$(A\times B)\wr G$ is a subgroup of $(A\wr G)\times (B\wr G)$. 
Now apply Lemma A and Proposition 2.3.\endprf

\proof{Proof of Main Theorem} Let $\Delta$ be a nontrivial
group with $\Gamma$ a strongly poly-surface normal subgroup of $\Delta$ of
finite index and $G=\Delta/\Gamma$. We will prove the theorem by 
induction on the rank of $\Gamma$. Note that $\Delta$ is a subgroup of
$\Gamma\wr G$. Hence it is enough to check that FIC is true for $\Gamma\wr
G$.

\noindent
{\bf Induction hypothesis $I(n)$}\quad For any strongly poly-surface group
$\Gamma$ of rank $\leq n$ and for any finite group $G$, FIC is true for
the wreath product $\Gamma\wr G$.

If the rank of $\Gamma$ is $\leq0$ then $\Gamma\wr G=G$ finite and
hence $I(0)$ holds.

Now assume $I(n-1)$. We will show that $I(n)$ holds.  

Let $\Gamma$ be a strongly poly-surface group of rank $\leq n$ and is
a normal subgroup of $\Delta$ with $G$ as the finite quotient group. So we
have a filtration by subgroups $$1=\Gamma_0\subset \Gamma_1\subset \cdots
\subset \Gamma_n=\Gamma$$ with all the requirements as in the definition
of strongly poly-surface group and there is the exact
sequence $$1\to \Gamma\to \Delta\to G\to 1$$
We have another exact sequence which is obtained after taking wreath 
product of the exact sequence $1\to \Gamma_1\to \Gamma\to
\Gamma/\Gamma_1\to 1$ with $G$. $$1\to
\Gamma_1^G\to \Gamma\wr G\to (\Gamma/\Gamma_1)\wr G\to 1$$
Let $p$ be the surjective homomorphism $\Gamma\wr G\to
(\Gamma/\Gamma_1)\wr G$. Note that $\Gamma/\Gamma_1$ is a strongly
poly-surface group of rank less or equal to $n-1$.

By induction hypothesis FIC is true for $(\Gamma/\Gamma_1)\wr G$. We
would like to apply Lemma C. Let $Z$ be a virtually cyclic subgroup of
$(\Gamma/\Gamma_1)\wr G$. Then there are two cases to consider.

\noindent
{\bf $Z$ is finite}\quad In this case we have $p^{-1}(Z) < \Gamma_1^G\wr
Z < \Gamma_1\wr (G\times Z)$. Since $\Gamma_1$ is a surface group, Lemma E
completes the proof in this case.

\noindent
{\bf $Z$ is infinite}\quad Let $Z_1=Z\cap (\Gamma/\Gamma_1)^G$. Then $Z_1$ is
an infinite cyclic normal subgroup of $Z$ of finite index. Let $Z_1$ be  
generated by $u$. We get $p^{-1}(Z) < p^{-1}(Z_1)\wr K$ where
$K$ is isomorphic to $Z/Z_1$. 

Also $$p^{-1}(Z_1)\wr K \ \simeq \ (\Gamma_1^G\rtimes \langle u\rangle
)\wr K < (\prod_{g\in G}(\Gamma_1\rtimes_{\alpha_g} \langle u\rangle ))\wr
K\eqno{\bf\label}$$
Now we describe the notations in the display (2.1). Let $t\in \Gamma^G$
which goes to $u$. Then $\alpha_g(\gamma)=t_g\gamma t_g^{-1}$ for all
$\gamma\in \Gamma_1$ and $t_g$ is the value of $t$ at $g$. By definition of
strongly poly-surface group each of these actions is induced by a
diffeomorphism of a surface $S$ whose fundamental group is isomorphic to
$\Gamma_1$.

Now there are two cases: (a) $\Gamma_1$ is finitely generated and (b)
$\Gamma_1$ is infinitely generated. 

\noindent
{\bf (a)}\qua Recall that if the fundamental group of a noncompact surface is
finitely generated then the surface is diffeomorphic to the interior of a
compact surface with boundary. Thus in this case the right hand side of
the display (2.1) is isomorphic to $(\prod_{g\in G}\pi_1(N^g))\wr K$ where
for each $g$, $N^g$ is a compact $3$-manifold fibering over the circle. If
$S$ is compact with nonempty boundary or is the interior of a compact
surface with nonempty boundary then $\partial N^g\neq\emptyset$ for all
$g$. In this situation use Corollary B to complete the proof of the
theorem. If $S$ is closed then so is $N^g$ for each $g$ and hence
Corollary 2.4 completes the proof if $S$ is not the Klein Bottle. If $S$
is the Klein bottle then the proof follows from the following lemma and
by noting that $\pi_1(S)$ has a finite index rank $2$ free abelian
subgroup. 

\proclaim{Lemma 2.5} Let $G_1,G_2,\cdots ,G_n$ be finitely presented
groups so that each $G_i$ contains a finitely generated free abelian
subgroup of finite index. For each $i$ let $f_i$ be an automorphism of
$G_i$. Let $G$ be a finite group. Then FIC is true for
the group $((G_1\rtimes_{f_1} \langle t \rangle )\times
(G_2\rtimes_{f_2} \langle t \rangle ) \times\cdots \times
(G_n\rtimes_{f_n} \langle t \rangle ))\wr G$.\endproc

\prf Recall that for groups $A,B$ and $G$, $(A\times B)\wr G$ is a
subgroup of $(A\wr G)\times (B\wr G)$. Also if FIC is true for two
groups then applying Lemma C twice and noting that FIC is true for
virtually poly-$\Bbb Z$ groups it follows that FIC is also true for the
product of the two groups. Thus it is enough to prove the Lemma for $n=1$. 

Note that by taking intersection of all conjugates of the
free abelian subgroup of $G_i$ we get a finitely generated free abelian
normal subgroup of $G_i$ with a finite quotient group, say $K_i$. 
Now since $K_i$ is a finite group and $G_i$ is finitely presented, there
are only finitely many homomorphism from $G_i$ onto $K_i$. Let $H_i$ be
the intersection of the kernels of these finitely many homomorphism. Then
$H_i$ is a finitely generated free abelian characteristic subgroup of
$G_i$ of finite index. Hence we have an exact sequence. $$1\to H_i\to
G_i\rtimes_{f_i} \langle t \rangle \to L_i\rtimes_{f_i} \langle t \rangle
\to 1$$ where $L_i \ \simeq \ G_i/H_i$. Taking wreath product with $G$ the
above exact sequence reduces to the following.
$$1\to H_i^G\to (G_i\rtimes_{f_i} \langle t \rangle)\wr G \to
(L_i\rtimes_{f_i} \langle t \rangle)\wr G \to 1$$  Note
that $(L_i\rtimes_{f_i} \langle t \rangle)\wr G$ is virtually
poly-${\Bbb Z}$ and $H_i^G$ is free abelian and hence FIC is true for
these two groups. Let $C$ be a virtually cyclic subgroup of
$(L_i\rtimes_{f_i} \langle t \rangle)\wr G$ then $p^{-1}(C)$ is easily
shown to be virtually poly-${\Bbb Z}$ and hence FIC is true for
$p^{-1}(C)$. Here $p$ denotes the last surjective homomorphism in the
above exact sequence. This completes the proof of the Lemma.\endprf

\noindent
{\bf (b)}\qua As $\Gamma_1$ is infinitely generated and free, by the
definition of strongly poly-surface group, $S$ has one end. Replacing
$\tilde F$ by $S$ in Lemma 2.2 we get $$(\prod_{g\in
G}(\Gamma_1\rtimes_{\alpha_g} \langle u\rangle ))\wr K <
\lim_{i \to \infty }((\prod_{g\in G}\pi_1(N^g_i))\wr K)$$
Now using Corollary B and  Theorem F we complete the proof of the
theorem.\endprf

\references
\Addresses\recd


\reflist
\key{FJ} {\bf F.T. Farrell} and {\bf L.E. Jones}, {\it Isomorphism
conjectures in
algebraic $K$-theory}, J. Amer. Math. Soc. 6 (1993) 249--297.

\key{Jo} {\bf L.E. Jones}, {\it A paper for F.T. Farrell on his 60'th
birthday}, preprint, SUNY at
Stony Brook, November 2002.

\key{R1}
{\bf Sayed K. Roushon},
{\it $K$-theory of virtually poly-surface groups}, Algebr. Geom. Topol.,
3 (2003), 103--116.

\key{R2}
{\bf Sayed K. Roushon},
{\it The Fibered isomorphism conjecture for complex manifolds},
March 2004, submitted for publication,
http://www.math.tifr.res.in/$\sim$ roushon/elliptic.html.

\key{R3}
{\bf Sayed K. Roushon},
{\it The Farrell-Jones isomorphism conjecture for $3$-manifold
groups}, May 2004, submitted for
publication, http://www.math.tifr.res.in/$\sim$
roushon/3-manifold-fic.html.

\endreflist\np

\def\startpage{1}\count0 1

\headline{\vbox to 0pt{\vskip -4.5mm\line{\small\phead\ifnum
\count0=\startpage ISSN 1472-2739 (on-line) 1472-2747 (printed)
\hfill {\pnum\folio}\else\ifodd\count0\def\\{ }%
\ifx\theshorttitle\relax\thetitle\else\theshorttitle\fi\hfill{\pnum\folio}
\else\def\\{ and }{\pnum\folio}\hfill\ifx\theshortauthors\relax\theauthors
\else\theshortauthors\fi\fi\fi}\vss}}
\footline{\vbox to 0pt{\vglue 0mm\line{\small\pfoot\ifnum\count0=\startpage
\copyright\ \gtp\hfill\else
\agt, Provisional Erratum\hfill\fi}\vss}}

\agt\hfill      
\hbox to 45truept{\vbox to 0pt{\vglue -13truept{\logomed A\kern -.37em{\logobig 
T}\kern -.38em G}\vss}\hss}
\break
{\small Provisional Erratum\nl 
Published: 12 August 2004}

\vglue 15pt
\cl{\large\bf Erratum to}
\smallskip
\cl{\large\bf `$K$-theory of virtually poly-surface
groups'}
\def\theaddress{School of Mathematics,
Tata Institute\\Homi Bhabha Road,
Mumbai 400005, India.}
\def\theemail{roushon@math.tifr.res.in}
\def\theurl{http://www.math.tifr.res.in/\char'176roushon/paper.html}
\medskip
{\leftskip25pt\rightskip25pt\small
{\bf Abstract}\qua In this note, we point out an error in the above
paper. We also refer to some papers where this error is corrected
partially and describe a positive approach to correct it
completely. 
\medskip
{\bf AMS classification}\qua 57N37, 19J10; 19D35

{\bf Keywords}\qua Strongly poly-surface
group, fibered isomorphism conjecture, $3$-manifold groups, pseudoisotopy
functor
\medskip}

In this note `FIC' stands for the Fibered Isomorphism
Conjecture of Farrell and Jones corresponding to the pseudoisotopy
functor (see [\ref{FJ}]).

In the proof of the main lemma of [\ref{R1}] we found some
filtration of the surface $\tilde F$ which is preserved by the
diffeomorphism $f$ and
used this filtration to find a filtration of the mapping torus $M_f$ of
$f$ by compact submanifolds with incompressible tori boundary. Recall
that $\tilde F$ was the covering of the surface $F$
corresponding to the commutator subgroup of $\pi_1(F)$ and $f:\tilde F\to
\tilde F$ was a lift of a diffeomorphism $g:F\to F$. Also recall that the
main lemma of [\ref{R1}] says that the FIC is
true for $\pi_1(M)$ where $M$ is the mapping torus of a
diffeomorphism of $F$.
The proof of the
existence of the above filtration
of $M_f$, we sketched in [\ref{R1}] is incorrect. In the proof of the main
lemma of [\ref{R2}] we show that some regular finite sheeted cover of
$M_f$
admits a
filtration of the required type provided $g$ satisfies certain conditions.
We called diffeomorphisms satisfying these conditions {\it special}
([\ref{R2}], section 1, definition). In fact if the diffeomorphism $g$ is
not
special then in general
such a filtration of a finite sheeted covering of $M_f$ may not exists.
Thus if we assume that $g$ is
special then a complete proof of the main lemma of [\ref{R1}]
is given in
[\ref{R2}].
For general $g$ we
prove the main lemma of [\ref{R1}]
in [\ref{R3}] assuming that the FIC is true for $B$-groups. By definition
a
$B$-group contains a finite index
subgroup isomorphic to the
fundamental group of a compact irreducible $3$-manifold with nonempty
incompressible boundary so that each boundary component is a
surface of genus $\geq 2$.

We have also proved in theorems 3.3 and 3.4 of [\ref{R3}] that the FIC
is true for a large class of $B$-groups. Here we mention that
the surjective part of the FIC for torsion free $A$-groups (see definition
3.1 of [\ref{R3}]) is already proved by L.E. Jones in [\ref{Jo}] and we
have
proved in proposition 9.3 of [\ref{R3}] that every $B$-group is an
$A$-group. For clarity we recall that an $A$-group contains a finite
index subgroup isomorphic to the fundamental group of a complete
nonpositively curved $A$-regular Riemannian manifold.

Also we should point out that in this situation the proof of
proposition 2.3 of [\ref{R1}] needs a slightly elaborate argument. We give
this proof in proposition 1.7 of [\ref{R2}]. Recall that proposition 2.3
says that the FIC is true for the fundamental group of a $3$-manifold
which has a finite sheeted cover fibering over the circle.

{\it Finally we record that at the time of writing this note, the
main lemma of [\ref{R1}] remains unproven in general and therefore the
proof of any
result where the main lemma is used, for general monodromy,
should be given an alternate argument. Also the proof of the main lemma
given in [\ref{R1}] is
withdrawn.}

\references\Addresses\recd

\bye